\overfullrule=0pt   \magnification=\magstep1  \input amssym.def
\font\huge=cmr10 scaled \magstep2
\def\De{\Delta}  \def\de{\delta}    \def\N{{\Bbb N}}
\def\Om{\Omega}  \def\oOm{{\overline{\Om}}}  \def\si{\sigma}
\def\Co{{{\cal P}_{>\ge}}}  \def\s{{\frak S}}
\def\qed{\quad{\vbox{\hrule\hbox{\vrule height5pt{$\,\,\,$}\vrule
height5pt}\hrule}}}
\def\Ct{{{\cal P}_{\le<}}}
\def\Ch{{{\cal P}_{><}}}
\def\Cf{{{\cal P}_{\le\ge}}}
\def\boxit#1{\vbox{\hrule\hbox{\vrule{#1}\vrule}\hrule}}
\def\splus{\,\,{\boxit{+}}\,\,}
\def\ssplus{\,{\boxit{$_{+}$}}\,}
{\nopagenumbers
\rightline{{June 1999}}\bigskip\bigskip
\centerline{{\huge The cyclic structure of unimodal permutations}}\bigskip
\bigskip
\centerline{{ Terry Gannon}}\bigskip
\centerline{{\it Department of Mathematical Sciences,}}
\centerline{{\it University of Alberta}}
\centerline{{\it Edmonton, Canada, T6G 2G1}}\bigskip
\centerline{{e-mail: tgannon@math.ualberta.ca}}\bigskip\bigskip

\noindent{{\bf Abstract}}\bigskip
Unimodal (i.e.\ single-humped) permutations may be decomposed into a
product of disjoint cycles.
Some enumerative results concerning their cyclic structure --- e.g.\ ${2\over 3}$
of them contain fixed points --- are given. We also obtain in effect a
kind of combinatorial universality for continuous unimodal maps, by
severely constraining the possible ways periodic orbits of any such
 map can nestle together. But our main observation (and
tool) is the existence of a natural
noncommutative monoidal structure on this class of permutations which respects
their cyclic structure. This monoidal structure is a little mysterious, and
can perhaps be understood by broadening the context, e.g.\
 by looking for similar structure in other classes
of `pattern-avoiding' permutations.\vfill\eject}\pageno=1

\noindent{{\bf 1. Introduction}}\bigskip

Let $\Delta(n)$ denote the set of all 
 unimodal permutations $\delta$ of $I_n:=\{1,2,\ldots,n\}$.
That is, for any such $\delta\in\De(n)$ there
exists an $m\in I_n$ satisfying\smallskip
\item\item{(i)} $a<b\le m\ \Rightarrow\ \de(a)<\de(b)$, and
\item\item{(ii)} $m\le a<b\ \Rightarrow\ \de(a)>\de(b)$.\smallskip

\noindent Of course, $m=\de^{-1}(n)$ is the maximum. Write $\De(\star)$
for $\cup_n \De(n)$, and $\s_n$ for the symmetric group.

$\De(\star)$ is the discrete analogue of the unimodal maps studied in
1-dimensional dynamical systems (see e.g.\ the classic [3]). For instance,
it has been observed that small populations have a tendency to grow,
and large ones decrease. The simplest model for this is a 
unimodal function. This was the motivation presented in [6] for
analysing the cyclic structure of $\delta\in\De(\star)$. More
generally, if $f$ is any continuous unimodal map and $J$ is any finite
set for which $f(J)=J$, then the restriction $f|_J$ is `topologically
conjugate' to a unique $\de\in\De(\star)$, called the `permutation
type' of $f|_J$ (explicitly, $\de=\Omega_J^{-1}\circ f\circ \Omega_J$
where $\Omega_J$ is an increasing bijection to be defined shortly). In
contrast to that of periodic orbits, the theory of
finite invariant sets $J$ of maps $f:I\rightarrow I$ is still largely
undeveloped --- a notable exception is [5] --- and this paper can be
regarded as a move in that direction for the special case of unimodal
maps. We return to this context in section 3.

Unimodal permutations also appear naturally in a second context.
We say that a permutation $\pi\in\s_n$ `contains' a pattern $\sigma=[\sigma 1,
\sigma 2,\ldots,\sigma k]\in
\s_k$  if there exist $k$ indices $1\le i_{\sigma 1}<i_{\sigma 2}<\cdots<
i_{\sigma k}\le n$ such that $\pi(i_1)<\pi(i_2)<\cdots<\pi(i_k)$, otherwise
we say $\pi$ `avoids' $\sigma$ [7]. Equivalently, $\sigma$ is contained in
$\pi$ iff the permutation matrix of $\sigma$ is a submatrix of the
permutation matrix of $\pi$.  For example, $[3,2,4,1]$ contains the
patterns $[2,1,3]$ (take the 3 indices $\{1,2,3\}$) and $[2,3,1]$
(take e.g.\  indices $\{1,3,4\}$), but avoids $[1,2,3]$ and $[3,1,2]$.
Write ${\frak S}_n(\sigma_1,\sigma_2,\ldots)$ for the set of all
$\pi\in{\frak S}_n$ avoiding all $\sigma_i$.
Questions involving pattern-avoidance arise for example in sorting problems
in computing science. A slightly more general notion: call a set $S\subseteq
\cup_n\s_n$ of permutations `closed' [1] if for any pattern
 $\si$ contained in any $\pi\in S$, we have
 $\si\in S$. Now, $\De(n)$ is precisely the set ${\frak
S}_n([213],[312])$ of all those permutations
which avoid both patterns $[2,1,3]$ and $[3,1,2]$, and
$\De(\star)$ is closed.

It is easy to calculate the cardinality of $\De(n)$. Note that
$$\|\{\de\in\De(n)\,|\,\de^{-1}(n)=m\}\|=\left({n-1\atop m-1}\right)\ ,$$
so we get $\|\De(n)\|=2^{n-1}$. For example, the four permutations in
$\De(3)$ are $[231]=(123)$, $[132]=(1)(23)$, $[321]=(13)(2)$ and $[123]=(1)(2)(3)$.

Considerably more difficult is the determination of the cardinality of
the {\it transitive} unimodal permutations --- the $n${\it -cycles}
---,  the set of which we will denote $\De_n$. For example, 
$\De_5$ consists of the cycles (12345), (13425) and (12435).
Weiss and Rogers [8], using methods related to [4], obtained 
$$\|\De_n\|={1\over n}\sum_{d|n\atop d\ {\rm odd}}\mu(d)\,2^{{n\over d}-1}\ ,
\eqno(1)$$
where $\mu$ is the M\"obius function. Thus about ${1\over n}$ of the
permutations in $\De(n)$ are transitive (of course, ${1\over n}$ is
also the corresponding fraction for ${\frak S}_n$). 
The formula in (1) appears in other contexts: for instance, it counts
the number of bifurcations of stable periodic orbits of the
quadratic family $x\mapsto x^2-a$. Write $\De_{\star}:=\cup_k\De_k$.

Let $J$ be any subset of ${\Bbb R}$ with cardinality $k$.
Define $\Om_J:I_k\rightarrow J$ to be the unique increasing bijection
from $I_k$ to $J$. We are mostly interested in $J\subset \N:=\{1,2,\ldots\}$,
in which case put $\overline{\Om}_J=\Om_{{\Bbb N}\setminus
J}:\N\rightarrow \N\setminus J$.

Any $\delta\in\De(n)$ decomposes uniquely of course into a
product of pairwise disjoint cycles. Each cycle will also be unimodal:
if $\de|_J$ is a cycle, then $\de_J:=\Om^{-1}_J\circ\de\circ\Om_J\in
\De_{\|J\|}$. We shall call $\de_J$ the {\it shape},  and $\|J\|$ the
{\it length}, of $\de|_J$. Cycles of length 1 are fixed-points.

In this paper we will investigate questions concerning the cyclic structure
of unimodal permutations --- see e.g.\ equations (2),(7) below. We shall
find, for example, that ${2\over 3}$ of all unimodal permutations have
fixed-points and ${2\over 5}$ have 2-cycles, compared with $1-e^{-1}$ and
$1-\sqrt{e}^{-1}$ of all permutations, respectively.
It will also be found that many combinatorial
properties of a cycle are independent of its shape.

Given any  $\de$, write $N_\de:\De_\star\rightarrow\{0,1,2,\ldots\}$
for the {\it cycle-counter}, where $N_\de(\de')$ equals the number of
cycles in $\de$ with shape $\de'$. For example $N_{(1)(26)(35)(4)}(12)=2$.
Note that for any $\delta\in\De(n)$, 
$n=\sum_k\sum_{\de'\in\De_k} kN_\de(\de')$.
                                   
A complementary question to computing $\|\De_k\|$ is, what is the number
${\cal N}(N)$ of           $\de$ with a given cycle-counter $N_\de=N$?
This question (actually a less fundamental one involving only the lengths
and not shapes of subcycles) was asked in [6].
The answer turns out to be simple:
$${\cal N}(N)=2^{\ell-1}\eqno(2)$$
where $\ell$ is the number of distinct $\delta'\in\De_\star$
with $N(\delta')\ne 0$. Such a simple answer should be hinting at some
deeper structure. Indeed, our proof of (2) will be constructive, in that we will
find an associative but noncommutative operation   `$\splus$' from
(a subset of) $\De(n)\times \De({n'})$ onto $\De({n+n'})$, obeying
$$N_{\de\ssplus\de'}=N_\de+N_{\de'}\ .\eqno(3)$$
Then (2) is essentially the statement that every $\de\in\De(\star)$ can
be built up uniquely from $\splus$.

Similar questions can be asked
for other pattern-avoiding classes of permutations. Standard practise
in combinatorics is to enumerate certain sets, and when two sets are
discovered to have the same cardinality, to try to establish an
explicit bijection between them. Not surprisingly, the focus here has been on
enumeration questions, although [1] has called for a structure theory
of `closed sets'. For instance, Knuth (1973) showed that for any
$\si\in{\frak S}_3$, the cardinality $\|{\frak S}_n(\si)\|$ equals the
$n$th Catalan number, while [7] showed $\|{\frak
S}_n([123],[132],[213])\|$ is the ($n+1)$-th Fibonacci number. Curiously, 
questions of cyclic structure appear to have been ignored, and yet
 pattern-avoiding classes
of permutations are precisely those classes for which
cyclic structure is natural to investigate --- that is, their
subcycles avoid those same patterns. We see in this paper that for at least
some such classes, e.g.\ the unimodal ones, we get interesting answers.
We briefly return to this in section 3.

Similarly, one could hope that other closed sets of permutations
would have a nice monoidal structure. Another obvious direction is to
try  to extend this theory to multi-modal permutations. Also
interesting should be (unimodal) {\it nonbijective} functions
$\gamma:I_n\rightarrow I_n$ --- these are considered e.g.\ in [5,6].

\bigskip\bigskip\noindent{\bf 2. The monoidal structure}\bigskip

Consider any 
$\de_1\in\De(k)$, $\de_2\in\De(\ell)$, and put $m_1=\de_1^{-1}(k)$ and
$m_2=\de_2^{-1}(\ell)$. Choose any $J\subset I_{k+\ell}$ with $\|J\|=k$
and write $\Om$ for $\Om_J$, $\overline{\Om}$ for $\overline{\Om}_{J}$.

By the `sum' $\de_1\oplus_J\de_2$ (or just $\de_1\oplus\de_2$ if $J$
is understood) we mean the permutation satisfying, for all $i\in I_{k+\ell}$,
$$(\de_1\oplus_J\de_2)(i)=\left\{\matrix{(\Om\circ \de_1\circ\Om^{-1})(i)&
{\rm if}\ i\in J\cr (\overline{\Om}\circ\de_2\circ\overline{\Om}^{-1})(i)&
{\rm if}\ i\not\in J\cr}\right.\ .\eqno(4)$$
That is, $\de_1$ and $\de_2$ are intertwined, with $\de_1$ placed at $J$.
For example,  $(13425)\oplus_{\{2,4,5,6,7\}}(123)= (138)(25647)$.

Obviously $N_{\de_1\oplus\de_2}=
N_{\de_1}+N_{\de_2}$. Of course, $m:=(\de_1\oplus\de_2)^{-1}
(k+\ell)\in\{\Om(m_1),\overline{\Om}(m_2)\}$. This operation `$\oplus$'
obeys a kind of commutativity and associativity: e.g.\ $\de_1\oplus_J\de_2
=\de_2\oplus_{I_{k+\ell}\setminus J}\de_1$. Of course any
 $\de\in\De(\star)$ can be written  as $\de=\de'\oplus_{J'}\de''\oplus_{J''}
 \cdots\oplus_{J^{(s-1)}}\de^{(s)}$ using obvious notation, where
$\de^{(i)}$ is the shape of the subcycle $\de|_{J^{(i)}}$ of $\de$.

Our immediate problem is, given $\de_1$ and $\de_2$, to find all $J$
such that $\de_1\oplus_J\de_2$ is unimodal. Our goal is Theorem 6.
Without loss of generality,
we will assume unless otherwise stated that $J$ obeys ${\Om}(m_1)
<\overline{\Om}(m_2)$. The following result
follows trivially from unimodality, and hints at how $\de_1$ and $\de_2$
must fit together.

\medskip\noindent{{\bf Lemma 1.}} {\it Assume $\de_1\oplus\de_2\in\De(k+\ell)$.
Then for any $a\in I_k$, $b\in I_\ell$,

\item{$\bullet$} $a> m_1$ and $b\ge m_2$ $\Rightarrow$ ${\Om}(a)> m$
and $\overline{\Om}(b)\ge m$;

\item{$\bullet$} $a\le m_1$ and $b< m_2$ $\Rightarrow$ ${\Om}(a)\le m$
and $\overline{\Om}(b)< m$;

\item{$\bullet$} $a> m_1$ and $b< m_2$ $\Rightarrow$ ${\Om}(a)>
\overline{\Om}(b)$;

\item{$\bullet$} $a\le m_1$ and} $b\ge m_2$ $\Rightarrow$ ${\Om}(a)<
\overline{\Om}(b)$.
\medskip

Simple as it is, this Lemma provides a major clue to the ideas which
follow. Indeed, partition the pairs $I_k\times I_\ell$ into 4 sets
$\Co$, $\Ct$, $\Ch$ and $\Cf$ defined in the obvious way (e.g.\ $\Co$
consists of all pairs $(a,b)$ where  $a> m_1$ and $b\ge m_2$). Our
approach is related to that of [4], except in how we must treat the
`turning points' $m_i$.

Specialise for now to 
$\de_1\in\De_k$, $\de_2\in\De_\ell$, and choose any $a\in I_k$,
$b\in I_\ell$. Define $a_i:=\de_1^i(a)$, $b_i:=\de_2^i(b)$ for
all $i\ge 0$. We shall consider the successive iterates $(a_0,b_0)$,
$(a_1,b_1)$, etc., up to the smallest $L$ for which $(a_L,b_L)\in\Ch\cup\Cf$ (such an
$L$ always exists by Lemma 2 below). From Lemma 1 we then know the relative
ordering of ${\Om}(a_L)$ and $\overline{\Om}(b_L)$. Now going backwards,
we can use unimodality (via Lemma 1) to determine the ordering
of ${\Om}(a_{L-1})$ and $\overline{\Om}(b_{L-1})$, then $\Om(a_{L-2})$
and $\oOm(b_{L-2})$, and ultimately $\Om(a_0)$ and $\oOm(b_0)$.
In this way we can (indirectly) find the unique set $J$, and hence the
unique unimodal sum $\de_1\oplus_J\de_2$. The next several lines,
culminating in Theorem 3, fill in this sketch.

Let ${S}=S(a,b)$ be the sequence whose $i$th term $S_i$ is 1, 2, 3,
or 4 if $(a_i,b_i)\in\Co$, $\Ct$, $\Ch$ or $\Cf$, respectively.
Call this sequence `finite of length $L\ge 0$' if $S_i\in\{1,2\}$ for all
$0\le i<L$ and $S_{L}\in\{3,4\}$.

\medskip\noindent{{\bf Lemma 2.}} {\it $S(a,b)$ is always finite.}\medskip

\noindent{{\bf Proof.}} Suppose for contradiction that for each $i$,
$a_i\le m_1$ iff $b_i<m_2$. Then clearly both $k,\ell>1$. Without loss of
generality take
$a_0=m_1$ (so $a_1=k$ and $b_0<m_2$), and let $n>0$ satisfy $b_n=\ell$. Then if
$a_n=k$,
we would have $a_{n-1}=m_1$ and $b_{n-1}=m_2$, contradicting our supposition.

Therefore our supposition yields both $a_1>a_n> m_1$ and $b_n>b_1\ge m_2$.
 Thus $a_{1+1}< a_{n+1}$ and $b_{1+1}>b_{n+1}$, which requires
$S_{1+1}=S_{n+1}$. This in turn implies $a_{1+2}<a_{n+2}$ iff $b_{1+2}>b_{n+2}$,
etc. Inductively, we get that $S$ has period $n-1$: $S_{1+i}=S_{n+i}$.
But then $S_{k\ell+n-1}=S_{n-1}=1$ contradicts $S_{k\ell}=S_0=2$.\qed\medskip

Let $S(a,b)$ be of length $L$, and let $M$ be the number of $\ell\le L$
such that $S_\ell$ is 1 or 3. Write $a\succ b$ if $M$ is odd, otherwise
$a\prec b$.

\medskip\noindent{\bf Theorem 3.} {\it Given $\de_1\in\De_k$ and $\de_2\in\De_\ell$,
there exists exactly one set $J$ with $\de_1\oplus_J\de_2\in\De(k+\ell)$,
satisfying $\Om(m_1)<\oOm(m_2)$.}\medskip

\noindent{\bf Proof.} By Lemma 1, if we have $\de_1\oplus_J\de_2\in\De(k+\ell)$,
then we immediately get:
$$\Om_J(a)<\oOm_J(b)\ \Longleftrightarrow\ a\prec b\ .\eqno(5)$$
There is at most one cardinality--$k$ set $J\subset I_{k+\ell}$ 
which  obeys (5). Conversely, by the
definition of $\prec$ and $\succ$, if we
find such a set $J$, then $\de_1\oplus_J\de_2$ will necessarily be
unimodal with $\Om_J(m_1)<\oOm_J(m_2)$. 
 Such a set will exist if, for any $a\in I_k,b\in I_\ell$, we have

\item\item{(i)} if $c\in I_k$, $c<a$ and $a\prec b$,
then $c\prec b$;

\item\item{(ii)} if $c\in I_\ell$, $c>b$ and $a\prec b$,
then $a\prec c$;

\item\item{(iii)} if $c\in I_k$, $c>a$ and $a\succ b$,
then $c\succ b$;

\item\item{(iv)} if $c\in I_\ell$, $c<b$ and $a\succ b$,
then $a\succ c$.

\noindent The proof of (i)-(iv) is a simplified version of the proof of Proposition
5 given below.\quad\qed\medskip

Of course, by `commutativity' of $\oplus$, we also get that there is
exactly one $J'$ for which $\de_1\oplus_{J'}\de_2\in\De(k+\ell)$ and
$\Om_{J'}(m_1)>\oOm_{J'}(m_2)$. Thus for transitive $\de_1\ne\de_2$,
there are precisely 2 distinct $\de\in\De(k+\ell)$ with cycles of shape
$\de_1,\de_2$; while for $\de_1=\de_2$ there will be exactly 1 such $\de$.
An important relation between the two $\de_1\oplus\de_2$ is given in Proposition
5(c) below.

For example, let $\de_1=(123)$, $\de_2=(13425)$. Then the sums of the
form $\de_1\oplus\de_2$, $\de_1\oplus\de_1$, $\de_2\oplus\de_2$ are
(138)(25647) (for $m_1\prec m_2$) and (148)(25637) (for $m_1\succ m_2$),
(136)(245), and $(1583\,10)(26749)$. These have $J=\{1,3,8\}$,
$\{1,4,8\}$, $\{ 1,3,6\}$ and $\{1,3,5,8,10\}$, respectively.

\medskip\noindent{{\bf Corollary 4.}} {\it Let $\de\in\De_n$, and write $m=\de^{-1}
(n)$. Let $J$ be the set in Theorem 3. Then $J$ 
contains exactly one element from $\{1,2\}$, one from $\{3,4\}$, etc.\ 
Moreover, $(\de\oplus
\de)^{-1}(2n)=2m$ if $n\equiv m$ (mod 2); otherwise $(\de\oplus
\de)^{-1}(2n)=2m-1$}.\medskip

Similar comments hold for    repeated sums $\de\oplus\cdots\oplus\de$.
To see the first assertion in Corollary 4, 
apply unimodality repeatedly to the inequalities $\ell\succ\ell+1$
and $\ell<\ell+1$ in order to produce a contradiction. The second assertion
follows by counting the number of times $\de^\ell(1)>m$ for $1< \ell<
n-2$ ($n-2$ is the length of $S(1,1)$), to determine whether or not
$1\succ 1$.

The following technical definition is crucial.\smallskip

\noindent{{\bf Definition.}} {\it Call $\de'\in\De_{k}$} {acute} {\it
if $n\equiv \de'{}^{-1}(k)$ (mod 2) (so the two
maxima of $\de'\oplus\de'$ run diagonally SW--NE like `$/$'), otherwise call
it} {grave}. {\it
Choose any $\de\in\De(n)$ and let $\de(J)=J$, and write $m(\de|_J)$ for the
maximum of $\de|_J$. By ${\cal A}(\de)$ we
mean the set of all acute $\de'\in\De_{\star}$ which are the shapes of
subcycles $\de|_J$ of $\de$; similarly, ${\cal G}(\de)$ will be the grave
shapes in $\de$. Write $\de'\in{\cal A}_>(\de)$ if $\de'\in{\cal A}(\de)$ and
there is some subset $J\subset I_n$ such that $\de|_J$ has shape $\de'$
and $m(\de|_J)>m(\de)$; define ${\cal A}_<(\de),{\cal G}_>(\de),{\cal G}_<
(\de)$ similarly.}\smallskip

For example, $(12\ldots k)$ is acute iff $k=1$.
For $\de=(1)(26)(35)(4)$, ${\cal A}_>(\de)={\cal
A}_<(\de)=\{(1)\}$, ${\cal G}_>(\de)=\{(12)\}$ and ${\cal G}_<(\de)=\emptyset$.

${\cal A}_<(\de)$ and ${\cal G}_>(\de)$ can be thought of as the  subcycles
of $\de$ of `positive type' [4] (or `orientation-preserving' subcycles).
Put another way, think of ${\cal
A}_<(\de), \ldots,{\cal G}_>(\de)$ as {\it multi-sets}, i.e.\ their
elements come with multiplicity. Then we will find that the
multiplicities in ${\cal A}_<$ and ${\cal G}_>$ can be arbitrarily
large, but those of ${\cal A}_>$ and ${\cal
G}_<$ can never exceed 1. Because of this, ${\cal A}_>$ and ${\cal
G}_<$  will play an important role in Theorem 6 below. 

We must generalise Theorem 3 by removing the transitivity requirement. This
is equivalent to considering multiple sums.

Select any $\de_i\in\De_{n_i}$, for $i=1,2,3$. We are
interested in constructing unimodal sums $\de=\de_1\oplus_{J_1}\de_2\oplus_{J_2}\de_3
$ of these three permutations which obey $\Om_1(m_1)<\Om_2(m_2)
<\Om_3(m_3)$, by applying the preceding
analysis to the partial sums $\de_{ij}:=\de_i\oplus\de_j$.
Here and elsewhere we write $\Om_i$ for $\Om_{J_i}$, and $m_i=\de_i^{-1}(n_i)$.
Define $\prec_{ij}, \succ_{ij}$ for $\de_{ij}$ as before.
We will require as usual that each $m_i\prec_{ij}m_j$.
Note that we have no hope to construct a unimodal
sum $\de_1\oplus\de_2\oplus\de_3$ with $\Om_1(m_1)<\Om_2(m_2)<\Om_3(m_3)$,
if both $n_1\succ_{12}n_2$ and $n_2\prec_{23}n_3$. We will find that this
is the only obstacle; to show that, we must establish the compatibility
of the orderings $\prec_{ij}$.

\medskip\noindent{{\bf Proposition 5.}} {\it Choose any $\de_i\in\De_{n_i}$,
 and let $\de_{ij}$ be as in the preceding paragraph.}

\item{(a)} {\it Assume that either $n_1\prec_{12}n_2$ or $n_2\succ_{23}n_3$.
Then for any $a\in I_{n_1}$, $b\in I_{n_2}$,
$c\in I_{n_3}$, we have both}\smallskip

\item\item{(i)} {\it $a\prec_{12}b$ and $b\prec_{23}c$ implies $a\prec_{13}c$;}

\item\item{(ii)} {\it $a\succ_{12}b$ and $b\succ_{23}c$ implies $a\succ_{13}c$.}

\medskip
\item{(b)} {\it There exist  sets $J_i$ such that $\Om_1(m_1)<\Om_2(m_2)
<\Om_3(m_3)$ and $\de_1\oplus_{J_1}\de_2\oplus_{J_2}\de_3$ is unimodal,
iff either $n_1\prec_{12} n_2$ or $n_{2}\succ_{23} n_3$. Moreover, when
such sets $J_i$ exist, they will be unique.}

\medskip\item{(c)} {\it Suppose $\de\in\De_k$, $\de'\in\De_\ell$, $\de\ne\de'$,
and let $\de\oplus_{J_A}\de'$ and $\de\oplus_{J_B}\de'$ be the two distinct
unimodal sums. Then $k\prec_A \ell$ iff $k\prec_B \ell$.

}
\medskip
\noindent{{\bf Proof of (a).}} Assume for contradiction that we have found
$a,b,c$ so that $a\prec_{12}b$ iff 

\noindent $b\prec_{23}c$ iff $a\succ_{13}c$. 
Write $a_\ell=\de_1^\ell(a)$, $b_\ell=\de_2^\ell(b)$, $c_\ell=\de_3^\ell(c)$,
 $m_{ij}=\de^{-1}_{ij}(n_i+n_j)$, and let $\Om_{ij}:I_{n_i}\rightarrow
I_{n_i+n_j}$, $\oOm_{ij}:I_{n_j}\rightarrow I_{n_i+n_j}$ be the  
increasing maps which build up $\de_{ij}$. 

Put $L_{13}$ for the length of the sequence $S_{13}(a,c)$ --- i.e.\
the smallest $0\le \ell<\infty$ such that $(a_\ell,c_\ell)\in \Ch\cup\Cf$.
Define $L'_{12}$ to be the smallest $0\le\ell\le
\infty$ such that either $\Om_{12}(a_\ell)>m_{12}>\oOm_{12}(b_\ell)$ or
$\Om_{12}(a_\ell)<m_{12}<\oOm_{12}(b_\ell)$. Define $L_{23}'$ similarly.
$L_{ij}'$ is the furthest point to which we can carry a recursive unimodality
argument for $\de_{ij}$.

Let $L={\rm min}\{L_{13},L'_{12},L'_{23}\}<\infty$.
For each $\ell<L$, we get by definition either:

\item{$\bullet$} $\Om_{12}(a_\ell)\le m_{12}$, $\oOm_{12}(b_\ell)\le m_{12}$,
$\Om_{23}(b_{\ell})\le m_{23}$,
$\oOm_{23}(c_\ell)\le m_{23}$, $\Om_{13}(a_\ell)\le m_{13}$, $\oOm_{13}(
c_\ell)\le m_{13}$; or

\item{$\bullet$} $\Om_{12}(a_\ell)\ge m_{12}$, $\oOm_{12}(b_\ell)\ge m_{12}$,
$\Om_{23}(b_{\ell})\ge m_{23}$,
$\oOm_{23}(c_\ell)\ge m_{23}$, $\Om_{13}(a_\ell)\ge m_{13}$, $\oOm_{13}(
c_\ell)\ge m_{13}$.

\noindent Therefore  unimodality repeatedly applied to ``$a_1\prec_{12}b_1$
iff $b_1\prec_{23}c_1$ iff $a_1\succ_{13}c_1$'' yields
$$a_L\prec_{12}b_L\ {\rm iff}\ b_L\prec_{23}c_L\ {\rm iff}\ a_L\succ_{13}c_L
\ .\eqno(6)$$

\noindent{{\bf Case i.}} $L=L_{12}'<L_{13}$.

\noindent $L<L_{13}$ means $a_L> m_1$ iff $c_L\ge m_3$. 
The $L=L_{12}'$ condition implies
$a_L> m_1$ iff $b_L< m_2$ iff $a_L\succ_{12}b_L$. Putting all this together
with (6) forces $b_L=m_2$ and
$n_2\prec_{23}n_3$, hence $n_1\prec_{12}n_2$, which contradicts
 $L=L_{12}'$.

\smallskip\noindent{{\bf Case ii.}} $L=L_{23}'<L_{13}$.

\noindent This is handled identically to {\bf Case i}.

\smallskip\noindent{{\bf Case iii.}} $L=L_{13}$.

\noindent $L=L_{13}$ means $a_L> m_1$ iff $c_L< m_3$, iff $a_L\succ_{13}
c_L$ iff  $b_L\prec_{23} c_L$ iff $a_L\prec_{12}b_L$. This forces
$a_L>m_1$, $c_L<m_3$,
$b_L= m_2$, $b_L\prec_{23} c_L$ and $a_L\prec_{12}b_L$, and hence both
$n_1\succ_{12}n_2$ and  $n_2\prec_{23}n_3$, contrary to hypothesis.

\smallskip\noindent {\bf Proof of (b)}.  Immediate from (a).

\smallskip\noindent {\bf Proof of (c)}.  Without loss of generality take
$k\ge \ell$, and suppose for contradiction $k\prec_A\ell$ but $k\succ_B\ell$.
Then by part (b), there exist sets $J_i$ such that $\gamma:=\de\oplus_{J_1}\de'
\oplus_{J_2}\de$ is unimodal and $\Om_1(m)<\Om_2(m')<\Om_3(m)$.
Write $a_\ell=\gamma^\ell(\Om_1(m))$, $b_\ell=\gamma^\ell(\Om_2(m'))$,
$c_\ell=\gamma^\ell(\Om_3(m))$.
The result follows from Corollary 4 and repeated unimodality: for each
$\ell$  we get either
$a_\ell<b_\ell<c_\ell$ or $a_\ell>b_\ell>c_\ell$, hence $k=\ell$
and $\de=\de'$.\qed\medskip

For any $\de\in\De_k,\de'\in\De_{\ell}$, write $\de\triangleleft\de'$ if
$\de\ne\de'$ and $k\prec \ell$ in $\de\oplus\de'$. Proposition 5 tells us
that this gives us a total-ordering on
$\De_{\star}$. The 1-cycle $(1)$ is the minimal element,
(12) is the second smallest, and there is no maximal element: in fact
$\de\triangleleft (12\ldots n)$ for any
$\de\in\De_k$, $\de\ne (12\ldots n)$,  with $\de^{-1}(k)<n$.
In fact this is precisely the ordering on $\De_\star$ discussed by
Metropolis-Stein-Stein (1973), and extended into a refinement of the
Sarkovskii ordering $3>_s5>_s\cdots>_s8>_s4>_s2>_s1$ of ${\Bbb N}$, by Baldwin et
al (see [2,5] and references therein). In particular, $\de\triangleleft
\de'$ iff any continuous map $f:I\rightarrow I$ having a periodic
orbit with permutation type $\de'$ will necessarily have another with
type $\de$ (fix $I=[0,1]$, say). In [2] this is extended to arbitrary
(i.e.\  nonunimodal)
cycles, where the ordering (called `forcing') is partial, and in [5]
forcing is further extended to arbitrary maps $\gamma:I_n\rightarrow
I_n$, where it is no longer antisymmetric. In the unimodal case,
everything is simpler.  Write
$\delta_{(k)}$ for min$\,\De_k$; e.g.\ for odd $k$,
$\de_{(k)}=(1,{n+1\over 2},{n+1\over 2}+1,{n+1\over
2}-1,\ldots,{n+1\over 2}-{n-3\over 2},{n+1\over 2}+{n-1\over 2})$. 
Then $k<_s\ell$ iff $\de_{(k)}\triangleleft\de_{(\ell)}$. 
A theorem of Bernhardt (1987), or our Corollary 4, implies that a given $\de\in\De_\star$ has
an immediate predecessor $\de'$ (with respect to `$\triangleleft$')
iff $\de$ is the `double' ${\cal D}\de'$ of $\de'$ (see 
the end of this section). These comments on $\triangleleft$ are not
used in this paper.

We are now prepared for the general theorem on $\oplus$.

\medskip\noindent{{\bf Theorem 6.}} {\it Let $\de_i\in\De(n_i)$, $i=1,2$.
Define $m_i=\de_i^{-1}(n_i)$, $J_i=\{1,\de_i(1),\de_i^2(1),\ldots\}$, and
$\widehat{\de}_i=\Om_i^{-1}\circ\de_i\circ\Om_i$. Then:}

\item{(i)} {\it if either ${\cal A}_>(\de_1)\cap{\cal A}_>(\de_2)$ or ${\cal G}_<
(\de_1)\cap{\cal G}_<(\de_2)$ are nonempty, then there are no unimodal sums
of the form $\de_1\oplus\de_2$;}

\item{(ii)} {\it if instead $\widehat{\de}_1\in{\cal G}_<(\de_2)$ or
$\widehat{\de}_2\in{\cal A}_>(\de_1)$,
then there is no unimodal sum $\de_1\oplus\de_2$ with $m_1\prec m_2$;}

\item{(iii)} {\it otherwise, there is exactly one unimodal sum
$\de_1\oplus \de_2$ with $m_1\prec m_2$.}\medskip

The proof follows from repeated application of Proposition 5(b). $\widehat{\de}_i$
is the shape of the subcycle of $\de_i$ containing the maximum.
Of course the analogous statements to those in Theorem 6(ii),(iii)
hold for  unimodal sums $\de_1\oplus\de_2$ with $m_1\succ m_2$.

 \medskip\noindent{{\bf Definition.}} {\it Given $\de_i\in\De(n_i)$, denote by}
  $\de_1\splus\de_2$ {\it the unique unimodal sum $\de_1\oplus_J\de_2$ obeying
  $m_1\prec m_2$ (when it exists).}\medskip

  We thus get a (partial) monoidal
structure on $\De(\star)$. It is associative but not commutative:

\medskip\noindent{{\bf Proposition 7.}} {\it Let $\de_i\in\De(\star)$.}

\item{(a)} {\it If both} $\de_1\splus\de_2$ {\it and} $\de_2\splus\de_1$
{\it exist, then they
will be equal iff $\widehat{\de}_1=\widehat{\de}_2$, using the
notation of  Theorem 6.}

\item{(b)} {\it If} $\de_1\splus(\de_2\splus\de_3)$ {\it exists, then so does}
$(\de_1\splus\de_2)\splus\de_3$ {\it and they are equal.}\medskip

We can extend the domain of `$\splus$' to all of
$\De(\star)\times\De(\star)$, in the following natural way. Define the {\it
double} ${\cal D}\de\in\De_{2k}$ of $\de\in\De_k$ to be
$$({\cal D}\de)(i)=\left\{\matrix{(\de\splus\de)(i)&{\rm if}\
i\not\in\{2m-1,2m\}\cr (\de\splus\de)(4m-1-i)&{\rm otherwise}\cr}\right.\ .$$
For example, ${\cal D}(12\ldots
k)=(1,3,\ldots,2k-1,2,4,\ldots,2k)$. It is a consequence of Corollary
4 that for any $\de\in\De_\star$, ${\cal D}\de$ is
the immediate successor of $\de$, and that $\de$ is acute iff
${\cal D}\de$ is grave.

Now, for any $\de\in\De(n)$ and $\de'\in\De_k$, define
$\de{\splus\!}_e\,\de'\in\De(n+k)$ by
$$\de{\splus\!}_e\,\de'=\left\{\matrix{\de\splus\de'&{\rm if}\
\de'\not\in{\cal A}_>(\de)\cr
({\oOm_J}^{-1}\circ\de\circ\oOm_J)\splus{\cal D}\de'&{\rm otherwise}\cr}
\right.\ ,$$
where the subcycle $\de|_J$ is the `obstacle' to forming
$\de\splus\de'$, i.e.\ the
subcycle of shape $\de'$ with $m(\de|_J)>m(\de)$. Conjugating
$\de$ by $\oOm_J$ squeezes out that subcycle. By associativity, this
defines the operator ${\splus\!}_e$ defined on all of
$\De(\star)\times\De(\star)$. `${\splus\!}_e$' is an associative extension
of `$\splus$': where $\splus$ exists, it equals ${\splus\!}_e$.
Although $\De(n){\splus\!}_e\,\De(n')=\De(n+n')$, equation
(3) will not always be satisfied. For example,
$(1)(26)(35)(4)\splus(13)(2)$ does not exist, but $(1)(26)(35)(4)
{\splus\!}_e\,(13)(2)=(1)(29)(38)(47)(56)$. We will use $\splus$ but
not ${\splus\!}_e$ in section 3.

\vfill\eject\noindent{{\bf 3. Discussion}}\bigskip

The monoidal structure `$\splus$' found in the previous section obeys
(3), by construction, and so of course is ideally suited for enumerations
involving cyclic structure in $\De(n)$. We give two examples.

When $\de_i\in\De(\star)$ are disjoint, i.e.\ don't have
any cycles with similar shapes, then both $\de_1\splus\de_2$ and $\de_2\splus
\de_1$ will be defined. Hence we get equation (2).

In comparison with (2), the number of permutations in the symmetric group $\s_n$ which
have precisely $n_k$ disjoint subcycles of length $k$ (so $n=\sum n_kk$)
is $n!/\prod_k k\cdot n_k!$.

Let $\de\in\De_k$, and call $\De_\de(n)$ the set of all $\de'\in\De(n)$
possessing a subcycle of shape $\de$: i.e.\ $N_{\de'}(\de)>0$. Then $\|\De_\de(n)
\|=(2^{n-k}-
2\cdot\|\De_{\de}(n-k)\|)+\|\De_{\de}(n-k)\|$, which can be solved to
yield
$$\|\De_\de(n)\|={1\over 2^k+1}(2^n-2^\ell(-1)^{[n/k]})\ ,\eqno(7)$$
where $0\le\ell<k$ obeys $n\equiv\ell$ (mod $k$), and $[x]$ is the
greatest integer not more than $x$. Thus about ${2\over 2^k+1}$ of
all unimodal permutations contain a given $\de\in\De_k$.

By comparison, we find that the number of permutations in $\s_n$ which
don't possess {\it any} $k$-cycles (when written as a disjoint product of
cycles) is precisely
$$n!\sum_{s=0}^{[n/k]}(-{1\over k})^s{1\over s!}\ ,$$
and thus their density converges to $e^{-{1\over k}}$.

Similar questions should be addressed for other pattern-avoiding sets
$\s_n(\si,\si',\ldots)$ of
permutations. For example, any $\pi\in{\frak S}_n([231],[312])$ is an
involution so is built from disjoint 1- and 2-cycles! Naturally, we
can't  expect all such sets to be equally interesting from this perspective
--- e.g.\ no permutations for $n>4$ can avoid both patterns $\{[123],
[321]\}$.
The choice $P=\{[123],[132]\}$ 
could be interesting to investigate from our point-of-view. Although there are
$2^{n-1}$ permutations which avoid $P$, as with $\De(n)$, there are 
 two 3-cycles which avoid $P$ (compared with $\|\De_3\|=1$), and
 both $(14)(23)$ and $(13)(24)$ have cycle structure $(12)\oplus(12)$
 (compared with only one unimodal sum $(12)\oplus(12)$).

It is important here to
 consider the following symmetry. It is known [7] that there
are 8 operations $\s_n\rightarrow \s_n$, forming
the dihedral group ${\frak D}_4$, that can be performed on our sets and which respect
questions of pattern-avoidance. In particular, we can hit any permutation on
the left or right with the involution $\iota=[n,n-1,\ldots,1]$, or we can
replace a permutation by its inverse: $\pi\mapsto \pi^{-1}$. For any choice of
operation $\alpha\in {\frak D}_4$, the set $\s_n(\alpha(\si),\alpha(\si'),\ldots)$
equals the set of all $\alpha(\pi)$ for $\pi\in\s_n(\si,\si',\ldots)$.
Half of these symmetries preserve in addition the cyclic structure:
 namely,  the four operations
$\pi\mapsto \pi,\pi^{-1},\iota\circ\pi\circ\iota,\iota\circ\pi^{-1}
 \circ \iota$, which together form a ${\Bbb Z}_2\times{\Bbb Z}_2$
 symmetry.
 
The unimodal permutations are precisely those which avoid both $\{[213]
,[312]\}$. Our ${\Bbb Z}_2\times{\Bbb Z}_2$ symmetry sends that to the
sets $\{[132],[231]\}$, $\{132],[312]\}$ and $\{[231],[213]\}$, so their
corresponding
pattern-avoiding sets will also possess a monoidal structure and satisfy
the same enumeration formulas (2),(7).

A question which seems to be relatively unexplored in the
1-dimensional dynamics literature is how distinct periodic orbits can
nestle together in a given continuous map.
This paper shows how severely constrained this is in the unimodal case.
 For instance, let
$c,c'\in{\rm Int}\,I$ be the turning points of unimodal maps
$f,f':I\rightarrow I$. Let ${\cal
O}_i=\{m_i,f(m_i),f^2(m_i),\ldots\}$, ${\cal O}_i' =\{m_i',f'(m_i'),
\ldots \}$ be sets of periodic orbits for $f$ and $f'$, where $m_i$ is
the maximum point of ${\cal O}_i$ (i.e.\ $f(m_i)={\rm max}\,{\cal
O}_i$), and similarly for $m_i'$. It is a consequence of our work that
{\it the finite bijections $f|_{\cup_i{\cal O}_i}$ and $f'|_{\cup_i
{\cal O}'_i}$ will have identical permutation type,} if for each $i$ $m_i$
and $m_i'$ have the same `itinerary' [4,3], i.e.\ (slightly more
strongly) if for each $i$,

\item{(i)} ${\cal O}_i$ and ${\cal O}_i'$ correspond to the same cycle
in $\De_\star$, and

\item{(ii)} either $m_i\le c$ and $m_i'\le c'$, or $m_i\ge c$ and
$m_i'\ge c'$.

For example, consider the nonconjugate maps $f(x)=0.939\sin\,\pi x$ and
$f'(x)=4x(1-x)$, and orbits ${\cal O}_1=\{0.5,.939,.179\}$, ${\cal
O}_2=\{.376,.869\}$, ${\cal O}_1'=\{.611,.950,.188\}$, and ${\cal
O}_2'=\{.345, .905\}$. Then the restrictions $f|_{{\cal O}_1\cup{\cal
O}_2}$ and $f'|_{{\cal O}_1'\cup{\cal O}_2'}$ are both conjugate to
the unimodal permutation (135)(24).

This observation can be regarded as a sort of combinatorial universality for
unimodal functions. Condition (ii) is related to the fact that
equation (2) involves a power of 2.

Consider now the logistic map $x\mapsto 4x(1-x)$.
All $\de\in\De_\star$ appear once or twice in it. 
A cycle will always appear there as a periodic
orbit of `negative type' or `orientation-reversing' [8]
 (i.e.\ with their maximum $<{1\over 2}$
for grave $\de$, and $>{1\over 2}$ for acute). Every $\de\in\De_n$ for
odd $n$ also appears as `positive type', but for even $n$ exactly
$\|\De_{n/2}\|$ (namely the doubles ${\cal D}(\De_{n/2})$)
 do not appear as positive type (this is a consequence of
[8]). For example, its fixed points are at $x=0$ (positive type) and
$x={3\over 4}$ (negative); its unique 2-cycle is $\{.345,.905\}$
(negative); and its 3-cycles are at $\{.188,.611,.950\}$ (positive)
and $\{.117,.413,.970\}$ (negative). That quadratic map thus implies the
existence (but not uniqueness) of many (but not all) sums $\de_1\oplus
\de_2\oplus \cdots$, and the ordering `$\triangleleft$' on $\De_\star$
can be read off from it --- e.g.\ since $0<.905<.950$, we have
$(1)\triangleleft (12)\triangleleft(123)$.

\bigskip \noindent{{\bf Acknowledgements}} \medskip

This research was supported
in part by NSERC. I benefitted from conversations with T.\ D.\ Rogers,
N.\ Lamoureux, V.\ Linek, and M.\ D.\ Atkinson.

\vfill\eject\noindent{{\bf References}}\bigskip

\item{[1]} M.\ D.\ Atkinson, Restricted permutations, Discrete Math.\
195 (1999) 27--38.

\item{[2]} S.\ Baldwin, Generalizations of a theorem of Sarkovskii on
orbits of continuous real-valued functions, Discrete Math.\ 67 (1987)
111--127. 

\item{[3]} P.\ Collet and J.-P.\ Eckmann, Iterated maps on the
interval as dynamical systems, Birkh\"auser, Boston, 1980.

\item{[4]} J.\ Milnor and W.\ Thurston, On iterated maps of the
interval, Lecture Notes in Mathematics, vol.\ 1342, Springer, Berlin,
1988, pp.\ 465--563.

\item{[5]} M.\ Misiurewicz and Z.\ Nitecki, Combinatorial Patterns for
Maps of the Interval, Memoirs Amer.\ Math.\ Soc.\ 456 (1991).

\item{[6]} T.\ D.\ Rogers and J.\ R.\ Pounder, Cycle structure in
discrete-density models, Discrete Appl.\ Math.\ 4 (1982) 199--210.

\item{[7]} R.\ Simion and F.\ W.\ Schmidt, Restricted permutations,
{Europ.\ J.\ Combin.}\ {6} (1985) 383--406.

\item{[8]} A.\ Weiss and T.\ D.\ Rogers, The number of orientation-reversing
cycles in the quadratic map, {CMS Conference Proc.}, Vol.\ 8,
Amer.\ Math.\ Soc., Providence, 1987, pp.\ 703--711.

\end